\documentclass[12pt]{amsart}
\newcommand{\bC}{\hbox{\bf C}_1}
\newcommand{\BC}{\hbox{\bf C}}
\newcommand{\bCC}{\hbox{\bf C}_2}
\newcommand{\bt}{{\mathbf 2}}
\newcommand{\bz}{{\mathbf 0}}
\newcommand{\bo}{{\mathbf 1}}
\newcounter{thlem}
\newtheorem{lemma}[thlem]{Lemma}
\newtheorem{theorem}[thlem]{Theorem}
\newtheorem{coll}{Corollary}[thlem]
\newcommand{\clos}{\hbox{clos\,}}
\newcommand{\isom}{\approx}
\let\ol\overline
\let\\\newline
\let\sp\vee
\let\inn\wedge
\let\cal\mathcal
\newcommand{\CM}{{\cal M}}
\newcommand{\COO}{\hbox{{\rm C${}_1$O${}_2$}}}
\newcommand{\CO}{\hbox{{\rm CO}}}
\newcommand{\Up}{{\cal{UP}}}
\newcommand{\Lo}{{\cal{LO}}}
\newcommand{\ideal}[2]{\pmb{\bigl[#1\bigr)}_{#2}}
\newcommand{\filter}[2]{\pmb{\bigl(#1\bigr]}_{#2}}

\title{Topological Representations of Posets}
\author[R.Breslav]{R.Breslav${}^1$}
\author[A.Stavrova]{A.Stavrova${}^1$}
\author[R.R.Zapatrin]{R.R.Zapatrin${}^2$}
\date{}
\begin{document}
\begin{abstract}
In~\cite{RRZ} an arbitrary poset $P$ was proved to be isomorphic
to the collection of subsets of a space $\CM$ with two closures
$\bC$ and $\bCC$, which are closed in the first closure and open in
other~-- $\COO(\CM,\bC,\bCC)$. As a space  for this representation
an algebraic dual space $P^*$ was used. Here we
extend the theory of algabraic duality for posets generalizing the notion
of an ideal.  This approach yields a sufficient condition for the
collection of $\COO$-subsets of $A\subset P^*$ (with respect to induced
closures) to be isomorphic to $P$. Applying this result to certain
classes of
posets we prove some representation theorems and get a topological
characterization of orthocomplementations.
\end{abstract}
\maketitle
\footnotetext[1]{The Centre of Mathematics, Nevsky pr., 39, 191011, St.Petersburg, Russia}
\footnotetext[2]{Department of Mathematics, SPb UEF, Griboyedova
30--32, 191023, St.Petersburg, Russia}
\section{Introduction}

Since Stone  introduced the topological representation of Boolean
algebras~\cite{Stone} there was a lot of attempts to generalize this
result: the Stone-like representations of orthopsets by Mayet and
Tkadlec~\cite{Mayet,Tkadlec},
different topological representations of distributive~\cite{Priestley,Rieger}
and arbitrary (by Hartonas, Dunn and Urquhart)~\cite{Hartonas-Dunn,Urq}
lattices. We follow the
construction introduced in~\cite{RRZ} where algebraic dual space
$P^*$ is endowed with two closures $\bC$ and $\bCC$ in such a way that
the collection of all subsets of $P^*$ which are closed in $\bC$ and
open in $\bCC$ ordered by set inclusion (we denote this collection
by $\COO(P^*,\bC,\bCC)$) is isomorphic to the initial poset $P$:

\begin{equation}\label{errz01}
\COO(A,\bC,\bCC)\isom P
\end{equation}

The representation (\ref{errz01}) of $P$ works for arbitrary poset
$P$. However, for particular classes of posets the `universal set'
$P^*$ can be contracted to a smaller one $A\subseteq P^*$ with the
closures $\bC$, $\bCC$ induced from $P^*$. In this paper we show
that the representations of specific classes of posets mentioned
above all have the form

$$
\COO(A,\bC,\bCC)\isom P
$$

\noindent and differ only by the choice of $A\subseteq P^*$.

\subsection{Spaces with two closures}

Mapping $\BC:\exp(\CM)\to\exp(\CM)$ we call {\em closure} if
\begin{enumerate}
\item $A\subset\BC(A)$;
\item $\BC(\BC(A))=\BC(A)$;
\item if $A\subset B$ then $\BC(A)\subset\BC(B)$.
\end{enumerate}

A set $A\subset\CM$ is {\em closed} (or {\em $\BC$-closed})
if $A=\BC(A)$, $A$ is {\em open}
if $\ol{A}=\CM\setminus A$ is closed and {\em clopen}
if it is both closed and open. Note, that any
intersection of closed sets is closed, and $\BC(A)$ is the
intersection of all closed sets which contain $A$. ${\cal
K}\subset\exp(\CM)$ is called the {\em base} of closure~$\BC$
($\BC=\clos(\cal K)$) if any closed set is an intersection of
elements of~$\cal K$.

The closure $\BC$ is {\em exact} if $\BC(\emptyset)=\emptyset$,
and {\em topological} if $\BC(A\cup B)=\BC(A)\cup\BC(B)$.
Note, that exact topological closure defines topology on $\CM$.
For a closure $\BC$ on $\CM$ define $\CO(\CM,\BC)$ to be the collection of
all clopen subsets of $\CM$. Obviously $\CO(\CM,\BC)$ ordered by set
inclusion is a bounded orthoposet. It was shown by Mayet and Tkadlec~\cite{Mayet,Tkadlec}, that
for an arbitrary bounded orthoposet $P$ there is a space $\CM$ with closure
$\BC$ such that $P\isom\CO(\CM,\BC)$.

If we define two closures on $\bC$ and $\bCC$ $\CM$, then by
$\COO(\CM,\bC,\bCC)$ we denote the collection of all subsets of $\CM$ which
are both $\bC$-closed and $\bCC$-open, ordered by set inclusion.
We can say nothing about the structure of $\COO(\CM,\bC,\bCC)$
except it is a poset, moreover, as it was shown in~\cite{RRZ} for
an arbitrary poset $P$ one can build a space with two closures such
that $P\isom\COO(\CM,\bC,\bCC)$.

\subsection{Algebraic duality for posets}\label{s12}

For a poset $P$ its {\em algebraic dual space} $P^*$ is the set
of all isotone mappings from $P$ to poset $\bt=\{0,1\}$ with $0<1$.
Here we develop the techniques needed to build the representation.

Consider $A\subset P^*$. A set $I$ we call an {\em ideal}
(with respect to $A$ or {\it $A$-ideal}) if $I$ is an intersection
of kernels of some mappings $x\in A$  (i.e. $I=\bigcap x^{-1}(0)$). Dually,
the intersection of co-kernels $F$ we call a {\em filter} ($F=\bigcap
x^{-1}(1)$).

For $B\subset P^*$ we define an ideal $I(B)$ (filter $F(B)$) to be the
intersection of kernels (co-kernels, respectively) of $x\in B$.

For $Q\subset P$ define an ideal $\ideal{Q}{A}$
(resp., filter $\filter{Q}{A}$)~-- the intersection of ideals
(resp., filters) containing $Q$.

Note that ideals with respect to $P^*$ coincide with order ideals ($I$ is
an order ideal if $q\in I$ and $p\leq q$ implies $p\in I$). In general
$A$-ideals are always order ideals, but the converse is not always
true.

We say that $A\subset P^*$ is {\em full} if for all $p\not\leq q$ there
exists $x\in A$ such that $x(p)=1$, $x(q)=0$.

$A\subset P^*$ is called {\em separating} if for any disjoint ideal
$I$ and filter $F$ there exists $x\in A$ such that $x\big|_I=0$ and
$x\big|_F=1$.

In some cases discussed in section~\ref{sec:com} $\filter{p}{A}$ and
$\filter{p}{A}$
coincide with lower and upper cones of $p$ respectively. Due to the
following obvious lemma the separating set is full in this case.

\begin{lemma}\label{lem:sep=>full}
Let $A$ be a separating subset of $P^*$ and
$\ideal{p}{A}\cap\filter{q}{A}=\emptyset$ for
all $q\not\leq p\in P$. Then $A$ is full.
\end{lemma}

\section{Topological representation: the general case\label{sec:gen}}

Define two closures on $P^*$. For $p\in P$ consider two subsets
of $P^*$:
$$
\Up(p)=\{x\mid x(p)=1\}\qquad\Lo(p)=\{x\mid x(p)=0\}.
$$
Then define closures $\bC$, $\bCC$ in the following way:
$$
\bC=\clos\{\Up(p)\}_{p\in P}\quad\mbox{and}\quad
\bCC=\clos\{\Lo(p)\}_{p\in P}.
$$
Note, that since $\ol{\Up(p)}=\Lo(p)$ all $\Up(p)$ are $\COO$-sets.

On $A\subset P^*$ consider closures ${\bC}_A, {\bCC}_A$ induced by
$\bC$ and $\bCC$ (i.e. ${\BC_i}_A(X)=\BC(X)\cap A$).
Let $\Up\!_A(p)=\Up(p)\cap A$ and $\Lo\!_A(p)=\Lo(p)\cap A$, then
$$
  {\bC}_A=\clos\{\Up\!_A(p)\}_{p\in P}\quad
  \mbox{and}\quad{\bCC}_A=\clos\{\Lo\!_A(p)\}_{p\in P}.
$$
We omit the index $A$ in ${\BC_i}_A$, $\Up_A$ etc.\
when it is clear which subspace is meant.

\medskip

The following equations show the relation between the closures
introduced on $A$ and $A$-ideals:
$$
{\bC}(X)=\bigcap_{p\in F(X)}\Up(p)\quad\hbox{and}\quad
\bCC(X)=\bigcap_{p\in I(X)}\Lo(p).
$$

\begin{theorem}\label{th:main}
Let $A\subset P^*$. Consider $\sigma:P\to\COO(A,{\bC}_A,{\bCC}_A)$ which maps $p$ to $\Up(p)$, then

\begin{itemize}
\item[$(1)$] $\sigma$ is isotone;
\item[$(2)$] if $A$ is full then $\sigma$ is injective;
\item[$(3)$] if $A$ is separating then $\sigma$ is surjective.
\end{itemize}
\end{theorem}
\begin{proof}
(1) Since $p\leq q$ implies $x(p)\leq x(q)$ for all $x\in P^*$
then $p\leq q$ implies $\Up(p)\subset\Up(q)$, so $\sigma$ is isotone.\\
(2) For $p\ne q$ either $p\not\leq q$ or $q\not\leq p$, so there
exists $x\in A: x(p)\ne x(q)$, then exactly one of $\Up(p)$, $\Up(q)$
contains $x$ and $\Up(p)\ne\Up(q)$.\\
(3) Let $B\in\COO(A,{\bC}_A,{\bCC}_A)$, then $B={\bC}_A(B)$ and
$\ol{B}={\bCC}_A(\ol{B})$. Consider $Q=I(\ol{B})\cap F(B)=I\cap F$. If
$Q=\emptyset$ there exists $x\in A: x\big|_I=0$ and $x\big|_F=1$, so
$x\in\Up(p)$ for all $p\in F$ and $x\in\Lo(q)$ for all $q\in I$. Thus
$x\in B$ and $x\in\ol{B}$ simultaniously, so $Q\ne\emptyset$. For $p\in Q$
we have $B\subset\Up(p)$, $\ol{B}\subset\Lo(p)=\ol{\Up(p)}$ and
$B=\Up(p)$.
\end{proof}

\begin{coll}
Let $A$ be a full and separating subspace of $P^*$, then
$P\isom\COO(A,{\bC}_A,{\bCC}_A)$.
\end{coll}

To get the topological representation of an arbitrary poset we prove

\begin{lemma}\label{lem:P*}
$P^*$ is full and separating.
\end{lemma}
\begin{proof}
For disjoint ideal $I$ and filter $F$, which are in this case order
ideal and filter, consider $x: x(p)=0$ for $p\in I$ and $x(p)=1$
otherwise. Obviously $x\in P^*$ and separates $I$ and $F$. Applying
lemma~\ref{lem:sep=>full} we see that $P^*$ is full.
\end{proof}

This leads us to the following theorem:

\begin{theorem}
Let $P$ be an arbitrary poset, then $P\isom\COO(P^*,\bC,\bCC)$.
\end{theorem}

\medskip

Due to the following lemma in the case of bounded poset $P$
subspaces $A$ of $P^*$ can be reduced:

\begin{lemma}
Let $P$ be a bounded poset, $A\subset P^*$ be full and separating,
then $A\setminus\{\bz,\bo\}$, where $\bz,\bo\in P^*$ are constant
mappings, is also full and separating.
\end{lemma}
\begin{proof}
Note that the ideals (filters) with respect to $A\setminus\{\bz,\bo\}$
coincide with the proper $A$-ideals ($A$-filters) and for disjoint nonempty
$I$ and $F$ the separating mapping $x\in A$ is not constant.
\end{proof}

\section{Topological representations: special cases\label{sec:com}}

We apply the results of previous section to some special classes of
posets.

\subsection{Orthoposets}

The bounded poset $P$ is called
an {\em orthoposet} if there exists an anti-isotone mapping
$(\cdot)':P\to P$ ({\em orthocomplementation}) such that $p=(p')'$,
$p\sp p'=1$ and $p\inn p'=0$. For an orthoposet define its {\em
orthodual} space $P^{*\prime}$ to be the set of all $x\in P^*$ such
that $x(p')=(x(p))'$.

\begin{lemma}
$P^{*\prime}$ is full and separating.
\end{lemma}
\begin{proof}
For disjoint ideal $I$ and filter $F$ consider $x: x(p)=0$ for $p\in I\cup
F'$, $x(p)=1$ for $p\in I'\cup F$, otherwise $x(p)=y(p)$ for some $y\in
P^{*\prime}$. Obviously $x\in P^{*\prime}$ and separates $I$ and $F$, so
$P^{*\prime}$ is separating. As $\ideal{p}{P^{*\prime}}$ is the lower cone of
$p$ for all $p\in P$ $P^{*\prime}$ is full according to lemma~\ref{lem:sep=>full}.
\end{proof}

Since $\Up\!_{P^{*\prime}}(p')=\Lo\!_{P^{*\prime}}(p)$, the bases
of closures $\bC$ and $\bCC$ coincide and $\bC=\bCC$. Denote

$$
\BC=\bC=\bCC
$$

\noindent Then $\COO$-sets are $\BC$-clopen. Applying
theorem~\ref{th:main} we have

\begin{theorem}\label{th:ortho}
Let $P$ be an orthoposet, then there exists a closure space $(\CM,\BC)$
such that $P\isom\CO(\CM,\BC)$.
\end{theorem}

\medskip

The representation obtained in previous theorem coicides with that
described by Mayet~\cite{Mayet} and Tkadlec~\cite{Tkadlec}.

Now we use the notion of full separating subspace to characterize
all orthocomplementations which can be introduced on a bounded
poset $P$. Any orthocomplementation $(\cdot)'$ defines a full
separating subspace of $P^*$ on which the closures $\bC$ and $\bCC$
coincide. Let $\cal S$ be the collection of full separating
subspaces of $P^*$ where $\bC=\bCC$. Consider $A\in\cal S$ then the
set complementation on $\COO(A,\bC,\bCC)\isom P$ is an
orthocomplementation, so with every $A\in\cal S$ we can associate
an orthocomplementation $(\cdot)'{}^A$ on $P$.

\begin{theorem}
All orthocomplementations on $P$ are in one-to-one correspondence with
maximal (with respect to set inclusion) elements of~$\cal S$.
\end{theorem}
\begin{proof}
For $A\in\cal S$ all $x\in A$ preserves $(\cdot)'{}^A$ because $x(p)=x(p'{}^A)=1$
implies $x\in\Up(p)$ and $x\in\Up(p'{}^A)=\ol{\Up(p)}$
(the similar contradiction holds for $x(p)=0$). It means that
$A\subset P^{*\prime A}$, so all maximal elements of $\cal S$ are of
the form $P^{*\prime A}$. Thus any orthodual space $P^{*\prime}$ is a
subspace of $P^{*\prime A}$ for some $A$. Obviously, orthocomplementation
associated with $P^{*\prime A}$ is $(\cdot)'{}^A$ and the one associated with
$P^{*\prime}$ is $(\cdot)'$. Since $P^{*\prime}\subset P^{*\prime A}$ and
orthocomplementations are induced by set complementation we get that
$(\cdot)'=(\cdot)'{}^A$ and $P^{*\prime}=P^{*\prime A}$, so all orthodual spaces, defined
by different orthocomplementations on $P$, are maximal in $\cal S$.
\end{proof}

\subsection{Distributive lattices}

According to the Stone representation theorem any Boolean algebra
is isomorphic to the collection of all clopen sets in some
topological space. Since Boolean algebra is an orthocomplemented
distributive lattice one can expect distributive lattice to be
represented as the collection of $\COO$-sets of some space with two
topological closures. We are going to construct such a
representation which follows from theorem~\ref{th:main} and is
different from Priestley~\cite{Priestley} and Rieger~\cite{Rieger}.

For a lattice $L$ let $L^{*\sp\inn}\subset L^*$ be the set of all
lattice morphisms (isotone mappings preserving lattice operations) from $L$
to $\bt$. Note that $L^{*\sp\inn}$-ideal is always lattice ideal
(an order ideal $I$ is called lattice ideal if $a,b\in I$ implies $a\sp b\in I$).

\begin{lemma}\label{lem:dist-full}
For any distributive lattice $L$ the ideals (filters)
with respect to $L^{*\sp\inn}$ coincide with the lattice ideals (filters).
Besides that, $L^{*\sp\inn}$ is full and separating.
\end{lemma}
\begin{proof}
First we prove that for disjoint lattice ideal $I$ and filter $F$
there exists $x\in L^{*\sp\inn}$ such that $x|_I=0; x|_F=1$ (it
means that $L^{*\sp\inn}$ separates lattice ideals). Suppose $I_0$
to be the maximal lattice ideal containing $I$ which is disjoint
with $F$.  The set-complement of $I_0$ is a filter~\cite{Gretzer},
thus the mapping $x$: $x|_{I_0}=0$, $x|_{L\setminus I_0}=1$
preserves $\sp$ and $\inn$.  For an arbitrary $p\in L$ the upper
cone of $p$ is a lattice filter. Then we get every lattice ideal
$I$ to be the intersection of kernels of all $x_p$, which separates
$I$ and the upper cone of $p$, over all $p\not\in I$, so I is an
ideal with respect to $L^{*\sp\inn}$ (recall the definition of
ideal in section \ref{s12}). Hence, the separating property for
$L^{*\sp\inn}$ is equivalent to the fact that $L^{*\sp\inn}$
separates lattice ideals, which was proved above.  $L^{*\sp\inn}$
is full by lemma~\ref{lem:sep=>full}.
\end{proof}

\begin{theorem}\label{th:dlat}
For any distributive lattice $L$ there exists a space with two
topological closures $(\CM,\bC,\bCC)$ such that $L\isom\COO(\CM,\bC,\bCC)$.
\end{theorem}
\begin{proof}
The only thing we need to prove is that the closures $\bC$, $\bCC$ induced
on $L^{*\sp\inn}$ are topological. Since elements of $L^{*\sp\inn}$
preserve both $\sp$ and $\inn$ we have $\Up(p\sp
q)=\Up(p)\cup\Up(q)$ and $\Lo(p\inn q)= \Lo(p)\cup\Lo(q)$, so the
bases of $\bC$ and $\bCC$ are closed under finite set union,
therefore the closures themselves are topological.
\end{proof}

\begin{coll}
A lattice $L$ is distributive iff $L^{*\sp\inn}$ is a full
separating subspace of $L^*$.
\end{coll}
\begin{proof}
This follows from lemma~\ref{lem:dist-full}, the fact that for any
lattice $L$ the closures induced on $L^{*\sp\inn}$ are topological,
and that for any space $\CM$ with two topological closures
$\COO(\CM,\bC,\bCC)$ is a distributive lattice.
\end{proof}

\subsection{Boolean algebras}

Here we present a proof of the Stone re\-p\-re\-sen\-ta\-tion
theorem:
\begin{theorem}[Stone] Any Boolean algebra $B$ is
isomorphic to the collection of all clopen subsets of a
topological space.
\end{theorem}
\begin{proof}
Since $B$ is a bounded distributive lattice,
$B^{*\sp\inn}\setminus\{\bz,\bo\}$ is full and separating.  Every
lattice morphism of Boolean algebras preserves
ortho\-comple\-mentation and, as in the case of orthoposets,
the topological closures $\bC$ and $\bCC$ do coincide.

Associating with every element of
$B^{*\sp\inn}\setminus\{\bz,\bo\}$ its kernel (that is a maximal
lattice ideal) one get the Stone space of Boolean algebra
originally described in~\cite{Stone}.
\end{proof}

\begin{coll}
Let $L$ be a distributive lattice, then $L$ is a Boolean algebra iff
closures $\bC$ and $\bCC$ coincide on
$L^{*\sp\inn}\setminus\{\bz,\bo\}$.
\end{coll}

\end{document}